\documentclass[12pt]{amsart}
\usepackage{amsthm}
\usepackage{url}
\newtheorem{lemma}{\bf Lemma}[section]
\newtheorem{proposition}[lemma]{\bf Proposition}

\newtheorem{definition}[lemma]{\bf Definition}
\gdef\cal{\mathcal}
\gdef\frak{\mathfrak}
\gdef\a{\aleph_1}
\begin{document}
\title{${\Omega \choose T}\neq {\Omega \choose \Gamma}$}
\author{Dominic van der Zypen}
\address{Federal Office of Borders and Customs Security,
Taubenstrasse 16, CH-3003 Bern, Switzerland}
\email{dominic.zypen@gmail.com}

\begin{abstract}
We show that the selection principles ${\Omega\choose T}$
and ${\Omega\choose\Gamma}$ are not equal constructing a topological
space $(X,\tau)$ that satisfies ${\Omega \choose T}$,
but not ${\Omega \choose \Gamma}$. This answers a
question from the year $2003$ in \cite{spm}.
\end{abstract}

\maketitle
\section{Introduction and definitions}
Throughout this note, let $(X,\tau)$ be a topological space.
\footnote{The definitions in this note also make sense in the
broader context of hypergraphs $H=(V,E)$, where
$V$ is any set and $E\subseteq\mathcal{P}(V)$.}
\begin{definition} \label{cov}
 We say that ${\cal U} \subseteq \tau$ is an {\em open cover}, or
  {\em cover} for short, if \begin{enumerate}
    \item $X\notin {\cal U}$,
    \item $\bigcup {\cal U} = X$.
  \end{enumerate}
\end{definition}
For $x\in X$ we set ${\cal U}_x = \{U \in {\cal U}:
x\in U\}$ and call ${\cal U}_x$ the
{\em star} of $x$ (with respect to ${\cal U}$).
\subsection{Thick covers} We call a cover ${\cal U}\subseteq \tau$
\begin{enumerate}
  \item {\em large} if ${\cal U}_x$ is infinite for every $x\in X$,
  \item an $\omega$-cover if every finite subset of $X$ is contained
    in some member of ${\cal U}$,
  \item a $t$-cover if ${\cal U}$ is large, and
    for all $x,y\in X$ at least one of the
    sets ${\cal U}_x\setminus {\cal U}_y$ and ${\cal U}_y\setminus
    {\cal U}_x$ is finite, and
  \item a $\gamma$-cover if ${\cal U}$ is infinite, and
    for every $x\in X$ the set ${\cal U}
    \setminus {\cal U}_x$ is finite.
\end{enumerate}
Let $\Lambda, \Omega, T, \Gamma$ denote the collections
of (open) large covers, $\omega$-covers, $t$-covers, and $\gamma$-covers,
respectively. An easy argument shows that every $\gamma$-cover is
large and therefore a $t$-cover, so $\Gamma \subseteq T$.

\subsection{The selection principle} If ${\mathfrak U}, {\mathfrak V}$ are
families of covers of $X$, then we define the property ${{\frak U}\choose
{\frak V}}$, read ``$\frak U$ choose $\frak V$'', as follows:

\begin{quote}
  ${{\frak U}\choose {\frak V}}\; :$ For each ${\cal U}\in {\frak U}$
  there is ${\cal V}\subseteq {\cal U}$ such that ${\cal V}\in {\frak V}$.
\end{quote}

\section{Construction of the example}
We consider the space $(\a, \tau)$ where $\a$ is the smallest
uncountable cardinal and $\tau$ is the collection of down-sets
in the cardinal $\a$, that is $\tau = \a\cup\{\a\}$.

\begin{proposition}\label{main}
  If ${\cal U}\subseteq (\tau \setminus \{\a\})$
  is a cover, then:
  \begin{enumerate}
    \item ${\cal U}$ is a large cover,
    \item ${\cal U}$ is an $\omega$-cover,
    \item ${\cal U}$ is a $t$-cover, but
    \item ${\cal U}$ is {\em never} a $\gamma$-cover.
  \end{enumerate}
\end{proposition}
{\it Proof.}

(1) Let $\alpha\in \a$. Suppose $\alpha$ is only covered by finitely
many members $U_1,\ldots U_n\in {\cal U}$. But then,
$\beta:=\bigcup{\cal U}\in \a$ is not covered by any member
of $\{U_1,\ldots, U_n\}$. We have $\beta > \alpha$, and since
${\cal U}$ is a covering of $\a$, there is $U^* \in {\cal U}$
covering $\beta$ and therefore $\alpha$. Clearly $$U^*\notin \{U_1,\ldots
, U_n\},$$ contradicting the assumption that the only members of
${\cal U}$ covering $\alpha$ are $U_1,\ldots, U_n$.

(2) Let $S\subseteq \a$ be finite, and consider
$$\beta = \bigcup S\in \a.$$ Then $\beta$ is contained in some
$U\in {\cal U}$, so $S \subseteq U$.

(3) Let $\alpha,\beta\in \a$. We may assume that $\alpha<\beta$.
Then ${\cal U}_\alpha\setminus
{\cal U}_\beta = \emptyset$, which is finite, so ${\cal U}$ is a
$t$-cover.

(4) Suppose that ${\cal U}\subseteq \tau\setminus\{\a\}$ is
a $\gamma$-cover. So every $\alpha\in\a$ is contained in all
but finitely many members of ${\cal U}$. Note that due to the
special nature of cardinals, we have ${\cal U}\subseteq \a$.
So ${\cal U}$ is a well-ordered set such that all members
only finitely many predecessors. This implies that ${\cal U}$
is either a finite or countable collection of members
of $\a$ such that $\bigcup{\cal U} = \a$.
This contradicts the fact that $\a$ is a {\em regular} cardinal.
$\Box$

By proposition \ref{main} (2) and (3), every $\omega$-cover in the space
$(\a,\tau)$ with $\tau = \a\cup\{\a\}$ is a $t$-cover, so the
property ${\Omega \choose T}$ is trivially true. On the
other hand, by proposition $\ref{main} (4)$, no cover is a $\gamma$-cover,
therefore property ${\Omega\choose \Gamma}$ is false.

So for the space $(\a, \a\cup\{\a\})$ we have ${\Omega\choose T} \neq
{\Omega\choose \Gamma}$.

As a final remark, the property ${\Omega\choose\Gamma}$
is the celebrated Gerlits-Nagy $\gamma$-property \cite{gamma}.
Since we have seen that $\Gamma \subseteq T$, property ${\Omega\choose
\Gamma}$ always implies ${\Omega \choose T}$. But the
converse is not true, as Proposition
\ref{main} shows that there is a space where ${\Omega\choose T}$ is
true but ${\Omega\choose \Gamma}$ is false.

{\bf Acknowledgement.} I am grateful to Boaz Tsaban for a fruitful
discussion via e-mail on the subject of this note.


\end{document}